\documentclass{ifacconf}

\usepackage{graphicx}      
\usepackage{natbib}        
\usepackage{amsmath}
\usepackage{amssymb}
\usepackage{dsfont}
\usepackage{xcolor}
\usepackage{caption}
\usepackage{subcaption}
\usepackage{algorithm}
\usepackage{algpseudocode}
\usepackage{url}

\newtheorem{problem}{Problem}
\newtheorem{lemma}{Lemma}
\newtheorem{theorem}{Theorem}
\newtheorem{remark}{Remark}
\newtheorem{assumption}{Assumption}
\newtheorem{example}{Example}
\newtheorem{definition}{Definition}

\newcommand{\rev}[1]{{\color{black} #1}}
\def\sk#1{\textcolor{blue}{#1}}

\begin{document}
\begin{frontmatter}

\title{Robust Data-Driven Receding Horizon Control} 
\author[First]{Jian Zheng} 
\author[Second]{Sahand Kiani} 
\author[First]{Mario Sznaier}
\author[Second]{Constantino Lagoa} 

\address[First]{ECE Dept., Northeastern University, 
   Boston, MA 02115 USA (e-mail: zheng.jian1@northeastern.edu, msznaier@coe.northeastern.edu)}
\address[Second]{EE Dept., Pennsylvania State University, State College, PA 16802 USA (e-mail: szk6437@psu.edu, cml18@psu.edu)}

\thanks{\copyright 2025 the authors. This work has been accepted to IFAC for publication under a Creative Commons Licence CC-BY-NC-ND.}
\thanks{J. Zheng and M. Sznaier were  partially supported by NSF grants CNS-2038493 and CMMI-2208182, AFOSR grant FA9550-19-1-0005, ONR grant N00014-21-1-2431  and DHS grant 22STESE0001-03-03.}





\begin{abstract}                
This paper presents a data-driven receding horizon control framework for discrete-time linear systems that guarantees robust performance in the presence of bounded disturbances. Unlike the majority of existing data-driven predictive control methods, which rely on Willem's fundamental lemma, the proposed method enforces set-membership constraints for data-driven control and utilizes execution data to iteratively refine a set of compatible systems online. Numerical results demonstrate that the proposed receding horizon framework achieves better contractivity for the unknown system compared with regular data-driven control approaches.
\end{abstract}


\end{frontmatter}
\section{Introduction}\label{sec:introduction}

Data-Driven Control (DDC) has gained significant attention in recent years due to the challenges of accurately modeling modern systems, as highlighted by \cite{hou2013model} and \cite{tang2022data}. DDC methods bypass explicit system identification, instead focusing on designing controllers that handle all systems consistent with observed data. Similarly, Model Predictive Control (MPC) has been widely adopted for its ability to handle constraints and optimize control inputs over a time horizon using system models as detailed in \cite{mayne2000constrained}. However, traditional MPC relies heavily on explicit system identification, which can introduce errors. Advances in data collection and computational capabilities have brought these two areas closer, enabling control design directly from measurement data. This paper builds on these developments, with DDC as the central focus and MPC playing a complementary role. Specifically, we aim to develop a robust Data-Driven Receding Horizon Control (DDRHC) method for discrete-time Linear Time-Invariant (LTI) systems. In the next subsection, we review relevant works on 
data-driven MPC.

\subsection{Related Works}
A study by \cite{berberich2020data} proposes a robust data-driven MPC method, based on Willems' Fundamental Lemma as detailed in \cite{willems2005note}, for LTI systems. The approach addresses bounded additive noise in output measurements by introducing a slack variable with cost regularization in the MPC problem. Nevertheless, the non-convexity introduced by the upper bound of the slack variable makes the problem computationally inefficient for practical implementation. Additionally, the new measurements are not utilized to refine the model, 
which may limit the performance. While most DDC studies are based on the Fundamental Lemma, in \cite{xie2024data}, the authors consider an unknown discrete LTI system affected by noise and design a general min-max MPC problem with input and state constraints. The considered min-max problem minimizes the worst case cost over the set of system matrices consistent with the measurement data. However, the authors restrict their approach to a linear state-feedback law and propose a receding horizon strategy by reformulating the data-driven min-max MPC problem as a semidefinite program.

Another approach to the data-driven MPC problem has been studied by \cite{attar2023data}, where they propose a method to compute robust backward reachable sets from noisy data for unknown constrained linear systems with bounded disturbances. 
Their approach constructs backward reachable sets using zonotopic inner approximations that align with measurements and state and input constraints, and designs a controller based on the reachable sets.
The study by \cite{li2023stochastic} introduces a data-driven receding-horizon control method for chance-constrained output tracking in unknown stochastic LTI systems influenced by both process and measurement noise. The approach constructs an auxiliary state model directly parameterized by input-output data, enabling the formulation of a stochastic control problem. This method incorporates chance constraints, aligning with the Stochastic MPC framework under investigation.

\cite{verhoek2021data} proposes a receding-horizon data-driven predictive control method for linear parameter-varying (LPV) systems, assuming the scheduling parameter is measurable and known over the prediction horizon. Using measured data and persistence of excitation, the method ensures reference tracking and constraint satisfaction with performance comparable to model-based control. However, the LPV system considered in their work is assumed to be noise-free.
\cite{schuurmans2021data} examines a receding horizon data-driven estimator for linear system networks, specifically focusing on partially observable Markov jump linear systems. Their proposed method uniquely identifies a subsequence of past mode transitions and removes the need for an offline exhaustive search over mode sequences to determine the observation window size.





\subsection{Contributions}

In this paper, unlike most DDC methods that rely on Willems' Fundamental Lemma—which often requires additional reformulation to achieve robust performance and may lead to non-convexities or inefficiencies in implementation—the proposed approach 
implements DDC by enforcing a set-membership constraint. Methods based on Willems' Fundamental Lemma typically produce a robust controller that is optimal for a single specific system but does not guarantee robust performance across all systems consistent with the measured data. In contrast, our method designs a robust controller that accounts for noise and ensures robust performance for all discrete-time linear systems compatible with the measured data. Furthermore, by utilizing a receding horizon framework, our approach updates the noise-affected data dictionary online at each time instant. Other methods may overlook state and input constraints, which limit their applicability to constrained systems, or may focus solely on linear state-feedback design laws without addressing designing an input in a general form. In contrast, our proposed method is applicable to all systems consistent with the measured data and dynamically integrates newly acquired data to refine the model in real-time. Additionally, our robust design formulation accounts for a general form of input while considering state and input constraints. The contributions of this work are as follows:

\begin{itemize}
    \item A novel robust data-driven receding horizon control framework that guarantees Uniformly Ultimately Bounded (UUB) stability and robust performance of all systems compatible with measured data under persistent noise, while state and input constraints are satisfied.
    \item A systematic approach for computing the largest robust controlled invariant set within state constraints, 
    ensuring feasibility.
    \item An efficient reformulation of the control problem using the extended Farkas' lemma is proposed, avoiding explicit vertex enumeration and thereby significantly reducing computational complexity. 
\end{itemize}

This paper is organized as follows: Section \ref{sec:introduction} outlines the motivation, challenges, and related works in DDC. Section \ref{sec:preliminaries} introduces key concepts like contractivity and boundedness. Section \ref{sec:robust_uub} formulates the robust UUB control problem under disturbances. Section \ref{sec:ddrhc} presents the data-driven receding horizon control framework, including invariant set handling. Section \ref{sec:main} proposes a reduced complexity framework of the optimization. Section \ref{sec:numerical_examples} provides numerical validation, and Section \ref{sec:conclusion} summarizes our contributions.

\section{Preliminaries}
\label{sec:preliminaries}

\subsection{Notation}
\begin{tabular}{p{0.15\columnwidth}p{0.75 \columnwidth}}
$\mathbb{R}, (\mathbb{R}^n)$ & Set of real numbers (of dimension $n$)\\
$x, X, \mathcal{X}$ & Vector (or scalar), matrix,  set of $x$\\
$\mathds{1},I$ & Vector of 1s, identity matrix\\
$\left\|x \right\|_\infty$ & $\ell_\infty$-norm of vector $x$\\
$\otimes$ & Kronecker product\\
\text{vec}$(X)$ & Vectorized matrix along columns \\ 
\text{int}$(\mathcal{X})$  & Interior  of set $\mathcal{X}$\\
$\partial \mathcal{X}$  & Boundary of set $\mathcal{X}$\\
\end{tabular}

\subsection{Minkowski Functionals and Ultimate Boundedness}
\begin{definition}
The Minkowski functional of a   set $\mathcal{P}$ is defined as
\[\psi_\mathcal{P}(x) = \text{inf}\left\{ r\in \mathbb{R}\colon r>0 \text{ and } x\in r\mathcal{P} \right\}.\]
\end{definition}
When the set $\mathcal{P}$ is convex, compact, symmetric
and contains the origin in its interior, $\psi(.)$ defines a norm in $\mathbb{R}^n$. In the sequel we will denote this norm by
$\|x\|_{\mathcal{P}}$.

Consider an uncertain discrete-time 
LTI system of the form:
\begin{equation}
    \label{eq:dynamics}
        x_{k+1} = A(w)x_k  + B(w)u_k +v_k,\; w \in \mathcal{W},\; v_k\in \mathcal{V},
\end{equation}
where $w$ and $v$ represent bounded uncertainty and process noise respectively, and $\mathcal{W},\mathcal{V}$ are compact, balanced polyhedra containing the origin in their interior. In the sequel, without loss of generality, we will assume that $\mathcal{V}$ has the representation $\mathcal{V}\doteq \{ v \colon \|Vv\|_\infty \leq 1\}$.

\begin{definition}[\hspace{1sp}\cite{BlanchiniUUB}] Given $0 \leq \lambda \leq 1$, a set $\mathcal{S}$ is (controlled) $\lambda$-contractive with respect to the trajectories of  \eqref{eq:dynamics} if and only if for all $x\in \mathcal{S}$ there exists a control $u =\phi(x)$ such that
$A(w)x  + B(w)u(x) +v \in \lambda \mathcal{S}$ for all $w \in \mathcal{W}, v \in \mathcal{V}$. In the case where $\lambda=1$, the set $S$ is said to be positively invariant. 
\end{definition}

\begin{lemma}[\hspace{1sp}\cite{BlanchiniUUB}] \label{lem:contractive} If a set $\mathcal{S}$ is  $\lambda$-contractive with respect to   \eqref{eq:dynamics}, then the set $\mu\mathcal{S}$ is also $\lambda$-contractive for all $\mu \geq 1$. Further, if $v_k \equiv 0$ then the result holds for all $\mu \geq 0$. 
\end{lemma}

\begin{definition}[\hspace{1sp}\cite{BlanchiniUUB}] The system \eqref{eq:dynamics} with a control of the form $u_k =\phi(x_k)$ is Uniformly Ultimately Bounded (UUB) in a set $\mathcal{S}$ if and only if
for every initial condition $x_o$, there exists $K(x_o)$ such that $x_k \in \mathcal{S}$ for all $k \geq K(x_o)$ and all $w \in \mathcal{W},v\in \mathcal{V}$.
\end{definition}
\begin{definition}
The system \eqref{eq:dynamics} with a control of the form $u=\phi(x)$  has a convergence rate $\lambda \leq 1$ to a set $S$
if and only if for all
$x\not \in \text{int}(\mathcal{S})$, 
$\|A(w)x  + B(w)\phi(x) +v\|_\mathcal{P}$
$ \leq \lambda\|x\|_\mathcal{P}$ for all 
$w \in \mathcal{W},v\in \mathcal{V}$.
\end{definition}

\begin{definition}
Given experimental data  $\mathcal{D}=\left\{x_k, u_k,\right.$ $\left .x_{k+1}\right\}_{k=0}^{T-1}$, consisting of $T$ state-input-next pairs, the consistency set is defined as
\begin{equation}\label{eq:cons}
\begin{aligned}
\mathcal{C}=&\left\{A,B\colon \|V(Ax_k+Bu_k - x_{k+1} \right )\|_\infty \leq 1, \\&  \left . k=0,\ldots,T-1 \right\}, 
\end{aligned}
\end{equation}
that is, the set of all systems consistent with the available information
(model structure, measured data and disturbance model). In the sequel, by a slight abuse of notation, we will denote this set as $\mathcal{C}(\mathcal{D})$.
\end{definition}

\subsection{Extended Farkas' Lemma}
The following variant of Farkas' Lemma plays a key role in this work to reduce the data-driven control design to a tractable convex optimization problem:

\begin{lemma}[\hspace{1sp}\cite{Henrion1999}]\label{lem:farkas} Consider polyhedrons
$\mathcal{P}_N \doteq \left \{ x \colon Nx \leq \nu \right \}$ and 
$\mathcal{P}_M\doteq \left \{ x \colon Mx \leq \mu \right \}$. Then $\mathcal{P}_N \subseteq \mathcal{P}_M$ if and only if there exists a matrix $Y$ with non-negative entries such that
\[
YN=M\; \text{ and } \; Y\nu\leq \mu.
\]
\end{lemma}

\section{Robust UUB Constrained Data-Driven Control}
\label{sec:robust_uub}


In this paper we consider uncertain systems of the form \eqref{eq:dynamics}, affected by persistent, bounded disturbances and subject to state and control constraints $x_k \in \mathcal{X}$, $u_k \in \mathcal{U}$. In the sequel, without loss of generality we will assume that the state and control constraint sets are given by:
\begin{equation}\label{eq:constraint}
\mathcal{X}\doteq \{ x \colon \|Fx\|_\infty \leq 1\} \;, \;
\mathcal{U}\doteq \{ u \colon \|Hu\|_\infty \leq 1\} \;
\end{equation}
Note that with this representation
$\psi_{\mathcal{X}}(x)=\|Fx\|_\infty$.

Due to the presence of persistent disturbances, it is only possible to guarantee that the closed-loop system can be rendered UUB in some set $\mathcal{S}$ that contains 
the smallest (controlled) origin reachable set from the disturbance.  This observation motivates the following prototype problem.

\begin{problem}\label{prob:proto} Given a $\lambda$-contractive set $\mathcal{S}$, find a control law $u = \phi(x)$ that  renders the closed-loop system UUB in $\mathcal{S}$ for all $w \in \mathcal{W}, v \in \mathcal{V}$ while respecting the state and control constraints.
\end{problem}

When the uncertain dynamics have 
\rev{a convex} description of the form $A(w)=\sum w_iA_i, \; B(w)=\sum w_i B_i$ where $A_i$,$B_i$ are known and $\mathcal{S}$ is polyhedral, the problem above can be solved using a polyhedral Lyapunov function induced by the Minkowski functional $\psi_\mathcal{P}$ of a suitable $\lambda$-contractive polyhedron $\mathcal{P} \subseteq \mathcal{S}$ (see \cite{BlanchiniUUB,blanchini2008set} for details). 

In the data-driven scenario considered here, \rev{however,} this approach cannot be applied directly, since $(A_i,B_i)$ are unknown and must be inferred from \rev{the observed data $\mathcal{D}$}.
Nevertheless, since $\mathcal{V}$ is a polyhedron, the consistency set $\mathcal{C}(\mathcal{D})$ is also a polyhedron \rev{whose vertices (i.e., $(A_i, B_i)$) can be easily computed}. Thus, in principle Problem \ref{prob:proto} can be solved by finding the vertices 
$(A_i,B_i)$ of $\mathcal{C}(\mathcal{D})$ and then proceeding as in \cite{blanchini2008set}. 
However, such an approach \rev{does} 
not make efficient use of all data available, 
\rev{i.e., the training data and the data generated by the current executions}, unless the consistency set $\mathcal{C}(\mathcal{D})$ is updated as new data is generated.  Furthermore, from a performance standpoint, it is desirable to minimize the size of the set $\mathcal{S}$ and achieve the fastest possible convergence rate compatible with the constraints. 
\rev{This motivates the following receding horizon problem.}

\begin{problem}\label{prob:RH1}
    Given experimental training data  $\mathcal{D}_{train}=\left\{x^d_j, u^d_j, x^d_{j+1}\right\}_{j=1}^{N}$  and data generated during the present execution up to the current time $T$,  $\mathcal{D}_{exe,T}=\left\{x_k, u_k,\right .$ $\left .  x_{k+1}\right\}_{k=0}^{T-1}$, find \rev{the smallest set $\mathcal{S}$ and a corresponding controller $u(x_T)$ that (i) renders all systems consistent with the data $\mathcal{D}_T=\mathcal{D}_{train}\cup \mathcal{D}_{exe,T}$  UUB in $\mathcal{S}$}
    compatible with the  state and control constraints, and (ii) has the fastest (
    \rev{i.e.,} smallest $\lambda$) convergence rate to this set.
\end{problem}

\section{Data-Driven Receding Horizon Control}
\label{sec:ddrhc}

A difficulty in solving Problem \ref{prob:RH1} stems from the fact  that computing an approximation to the smallest controlled invariant (or controlled contractive) set \rev{$\mathcal{S}$} is nontrivial \rev{as shown in} \cite{Ozay2018}. In addition, optimizing the rate of convergence requires optimizing the worst case (over initial conditions and disturbances) of  $\prod_{k=0}^T \lambda_k$, where $\lambda_k$ denotes the contractivity rate at time $k$, and $T$ is the considered time horizon. To address these challenges, we adopt a greedy receding horizon approach. At each time instant $k$, we minimize the worst-case $\lambda_k$ across all systems consistent with the available data $\mathcal{D}_k$ up to time $k$, while satisfying the control and state constraints.

\subsection{\rev{The Case with Controlled Contractive $\mathcal{X}$}}
\label{sec:Xinvariant}

For simplicity, first consider the case where the \rev{entire} state constraint set $\mathcal{X}$ is 
controlled $\lambda$-contractive, 
for some 
$0 \leq \lambda \leq 1$. In this case, 
\rev{the following algorithm is used as a surrogate of Problem \ref{prob:RH1}}:

\begin{algorithm}
\caption{Conceptual Data-Driven Receding Horizon Control}\label{alg:conceptual}
\begin{algorithmic}[1]

\State \textbf{Initialize:} initial state $x_0$ and data $\mathcal{D}_0=\mathcal{D}_{train}$
\Repeat
\State
\textbf{Input:} current state $x_k$ and data $\mathcal{D}_k$
\State Solve 
\begin{equation} \label{eq:conceptual}
\begin{aligned}
    u_k & = \mathop{argmin}_{u \in \mathcal{U}, \lambda}\;  \lambda  \quad \text{subject to:}\\
    & 
    \psi_\mathcal{X}(Ax_k+Bu+v) \leq  \lambda \psi_\mathcal{X}(x_k),  \\
    &\forall v \in \mathcal{V}, \; (A,B)\in \mathcal{C}(\mathcal{D}_k)
\end{aligned}
\end{equation}
\State Apply the control $u_k$ to the unknown system
\State Measure the resulting state $x_{k+1}$
\State Update the data set: $\mathcal{D}_{k+1} \gets \mathcal{D}_k \cup (x_k,u_k,x_{k+1})$
\State $k \gets k+1$
\Until{a stopping criterion is met}
\end{algorithmic}
\end{algorithm}

\begin{lemma}\label{lemma:stable} Assume that $\mathcal{X}$ is controlled $\lambda_o$-contractive for some 
$0 <\lambda_o \leq 1$, for all systems in
$\mathcal{C}(\mathcal{D}_{train})$. Then
the control law generated by Algorithm 1 is admissible and such that:
\begin{enumerate}
    \item
    The origin is an exponentially stable point of the \rev{nominal} closed-loop system, with region of attraction $
\mathcal{X}$ and convergence rate $\lambda^*$, for some  $\lambda^*< \lambda_o$.
\item The closed-loop system is UUB in the set $
\tilde{\lambda}  \mathcal{X} \subseteq \lambda_o\mathcal{X}$, with $\tilde{\lambda} \leq \lambda_o$.
\end{enumerate}
   \end{lemma}
   \begin{pf} 
   To prove (1), set $v_k \equiv 0$ and denote by $x_j$ the vertices of $\mathcal{X}$, where we assume a finite number of extreme points. 
   Since $0 \in \text{int}(\mathcal{V})$, it follows that for each $j$ there exists $u(x_j)$  and $\lambda_j^* < \lambda_o$ such that for all $(A,B) \in \mathcal{C}(\mathcal{D}_{train})$, $\psi_{\mathcal{X}}(Ax_j+Bu(x_j))\leq \lambda_j^*$.  From convexity, it follows that for all $x\in \mathcal \partial X$, there exists $u(x)$ such that $\psi_{\mathcal{X}}(Ax+Bu(x)) \leq \lambda^* \doteq \max_j\lambda_j^* < \lambda_o $. Consider now the trajectory generated by Algorithm \ref{alg:conceptual}, starting from some $x_0 \in \mathcal{X}$. Since at time $k$, $\mathcal{C}(\mathcal{D}_{k}) \subseteq \mathcal{C}(\mathcal{D}_{train})$, from Lemma \ref{lem:contractive} it follows that the control action generated by \eqref{eq:conceptual} is such that
   $\psi_\mathcal{X}(x_{k}) \leq \lambda^* \psi_\mathcal{X}(x_{k-1}) \leq (\lambda^*)^{k} \psi_\mathcal{X}(x_0) \leq (\lambda^*)^k$. Property (2) follows directly from Lemma  \ref{lem:contractive}  and the fact that, at each step, $\mathcal{C}(\mathcal{D}_{k}) \subseteq \mathcal{C}(\mathcal{D}_{train})$. Hence, if $\mathcal{X}$ is controlled $\lambda_o$-contractive for all $(A,B) \in \mathcal{C}(\mathcal{D}_{train})$, it is controlled $\tilde{\lambda}$-contractive for all
    $(A,B) \in \mathcal{C}(\mathcal{D}_{k})$, for some  $\tilde{\lambda} \leq \lambda_o$.
   \end{pf}
   \begin{remark} Note that condition (2) above is the tightest worst-case UUB condition achievable, since it encompasses the case where the true system is the one corresponding to the worst contractivity rate.
   \end{remark}
\subsection{The General Case} \label{sec:general}
In the case where the state constraint set $\mathcal{X}$ is not controlled contractive, there is no control action $u \in \mathcal{U}$ that guarantees satisfaction of the constraints for all $x \in \mathcal {X}$ and $v \in \mathcal{V}$. Nevertheless, Algorithm \ref{alg:conceptual} can still be applied by restricting the initial condition to $\mathcal{X}_I$, the largest controlled invariant 
\rev{subset of  $\mathcal{X}$,} and replacing $\psi_{\mathcal{X}}(.)$ by $\psi_{\mathcal{X}_I}(.)$ in the optimization \eqref{eq:conceptual}.

The set $\mathcal{X}_I$ can be constructed by applying the algorithm proposed in  \cite{blanchini2008set} to the vertices $(A_i,B_i)$ of the consistency set $\mathcal{C}(\mathcal{D}_{train})$. For ease of reference, this algorithm is shown below: 
\begin{algorithm} 
\caption{Backward Construction of Invariant Sets}\label{alg:invariant}
\begin{algorithmic}[1]
\State \textbf{Input:} $\mathcal{X}(F,g) = \left\{x \colon Fx \leq g\right\},\mathcal{U}(H) = \left\{u \colon Hu \leq 1\right\}$
\State \textbf{Initialize:} $k = 0,\; F^{(k)} = F,\; g^{(k)} = g,\; \mathcal{X}^{(k)} = \mathcal{X}$
\Repeat
\State Compute $\tilde{g}^{(k)}_i = g^{(k)}_i - \max_{v \in \text{vert}(\mathcal{V})} F^{(k)}_i v$
\State $\mathcal{M}^{(k)} = \left\{(x, u) \colon F^{(k)}[A_i x + B_i u] \leq \tilde{g}^{(k)},\; Hu \leq 1\right\}$
\State Compute $\mathcal{R}^{(k)} = \left\{x \colon \exists u,\; (x, u) \in \mathcal{M}^{(k)}\right\}$
\State Update $\mathcal{X}^{(k+1)} = \mathcal{R}^{(k)} \cap \mathcal{X}^{(k)}$
\State $k \gets k + 1$
\Until{a stopping criterion is met}
\State \textbf{Output:} $\mathcal{X}^{(k)}$
\end{algorithmic}
\end{algorithm}

\section{A Reduced Complexity Formulation}\label{sec:main}

Let $A_i^{(k)},B_i^{(k)}, i=1, \ldots, N_v$ denote the vertices of the set $\mathcal{C}(\mathcal{D}_k)$.
Since all sets involved are polyhedral,  a necessary and sufficient condition for feasibility of \eqref{eq:conceptual} is the existence of a control $u \in \mathcal{U}$ such that 
\begin{equation}\label{eq:ns}
     \psi_\mathcal{X}(A_i^{(k)}x_k+B_i^{(k)} u+v) \leq  \lambda \psi_\mathcal{X}(x_k),\; \forall v \in \mathcal{V}, \; i=1,\ldots, N_v \\
\end{equation}
Thus, in principle the optimization \eqref{eq:conceptual} reduces to a linear program. However, this requires finding the vertices of the consistency set $\mathcal{C}(\mathcal{D}_k)$ at each step, which is a nontrivial task.
Further, it  can lead to a very large number of constraints, since, as shown in \cite{avis2018mplrs}, the number of vertices of the consistency set grows combinatorially in 
$n$ and $T$.

To circumvent this difficulty, we propose an alternative approach based on duality that does not require finding the vertices of $\mathcal{C}(\mathcal{D}_k)$. 
Using the Kronecker product, $\mathcal{C}(
\mathcal{D}_k)$ can be equivalently reformulated as a polytope in the space of system coefficients:
\begin{equation} \label{eq:p1}
    \mathcal{C}(\mathcal{D}_k) = \left\{ a,b:
\begin{bmatrix}
    P_1 & Q_1 \\ -P_1 & -Q_1
\end{bmatrix}
\begin{bmatrix}
    a \\ b
\end{bmatrix} \leq 
\begin{bmatrix}
     \mathds{1}+\xi \\  \mathds{1}-\xi
\end{bmatrix}
\right\},
\end{equation}
where $a = \text{vec}(A^T)$, $b = \text{vec}(B^T)$ and
\[
P_1 = 
\begin{bmatrix}
    V \otimes [x^d_1]^T\\
    \vdots\\
     V \otimes [x^d_{N-1}]^T \\ \vdots \\
    V \otimes x_{k-1}^T\\
\end{bmatrix}, 
Q_1 = 
\begin{bmatrix}
 V \otimes [u^d_1]^T\\
    \vdots\\
     V \otimes [u^d_{N-1}]^T \\ \vdots \\
    V \otimes u_{k-1}^T\\
\end{bmatrix}, 
\xi = 
\begin{bmatrix}
    Vx_2^d \\
    \vdots \\
    Vx_{N}^d \\
    \vdots \\
    Vx_k
\end{bmatrix}
\]

We make the following assumption on the amount of data required for our data-driven control design:
\begin{assumption}
    Sufficient data are collected such that the matrix $[P_1, Q_1]$ has full column rank.
\end{assumption}
This assumption ensures that the consistency set $\mathcal{C}(\mathcal{D}_k)$ is compact. Note that this assumption is also necessary for identification-based methods. Otherwise, the worst-case identification error of any interpolation algorithm becomes unbounded, thereby causing identification-based robust control to fail.
\begin{theorem} \label{thm:conceptual} The optimization \eqref{eq:conceptual} is equivalent to the following linear program:
\begin{equation} \label{eq:dual_form}
\begin{aligned}
    u_k = \mathop{argmin}_{u \in \mathcal{U},Y \geq 0, \lambda}\;  \lambda  \quad &\text{subject to:}\\
Y
\begin{bmatrix}
    P_1 & Q_1 & 0\\ -P_1 & -Q_1 & 0 \\
    0 & 0 & V \\
    0 & 0 & -V
\end{bmatrix} &= 
\begin{bmatrix}
    F(I \otimes x^T) &  F(I \otimes u^T) &F \\ - F(I \otimes x^T) & - F(I \otimes u^T)& -F
\end{bmatrix} \\
 Y \begin{bmatrix}
     \mathds{1}+\xi \\  \mathds{1}-\xi \\ \mathds{1} \\ \mathds{1}
\end{bmatrix}
&\leq  \lambda \psi(x) \mathds{1}
\end{aligned}
\end{equation}
\end{theorem}
\begin{pf}
For a given $0\leq \lambda \leq 1$, let 
\[\mathcal{P}_2(u,x,\lambda) \doteq \{A,B,v \colon \|F(Ax+Bu+v)\|_\infty \leq \lambda \psi(x) \}\]
Using the properties of the Kronecker product yields the equivalent representation
\begin{equation} \label{eq:p2}
\begin{aligned}
\mathcal{P}_2 = \left\{ \begin{bmatrix}
    a \\ b \\v
\end{bmatrix}: \right. &
\begin{bmatrix}
    F(I \otimes x^T) &  F(I \otimes u^T) &F \\ -F(I \otimes x^T) & - F(I \otimes u^T)& -F
\end{bmatrix}
\begin{bmatrix}
    a \\ b \\v
\end{bmatrix} \\
& \leq \lambda \psi(x) \mathds{1} \bigg \}.
\end{aligned}
\end{equation}

Consider now the polyhedron defined by the  consistency set:
\[
\begin{aligned}
\mathcal{P}_1 &\doteq \mathcal{D}_k \times \mathcal{V} =\\
& \left \{\begin{bmatrix}
    a \\ b \\v
\end{bmatrix} : \begin{bmatrix}
    P_1 & Q_1 & 0\\ -P_1 & -Q_1 & 0 \\
    0 & 0 & V \\
    0 & 0 & -V
\end{bmatrix}
\begin{bmatrix}
    a \\ b \\v
\end{bmatrix} \leq 
\begin{bmatrix}
     \mathds{1}+\xi \\  \mathds{1}-\xi \\ \mathds{1} \\ \mathds{1}
\end{bmatrix}
\right\}
\end{aligned}\]
It can be easily seen that the constraint in \eqref{eq:conceptual} holds for a given $\lambda$ if and only if $\mathcal{P}_1 \subseteq \mathcal{P}_2$.
Applying the extended Farkas' Lemma, this inclusion is equivalent to the existence of a multiplier matrix $Y$ with non-negative elements such that
\begin{equation}\label{eq:duality}
\begin{aligned}
Y
\begin{bmatrix}
    P_1 & Q_1 & 0\\ -P_1 & -Q_1 & 0 \\
    0 & 0 & V \\
    0 & 0 & -V
\end{bmatrix} &= 
\begin{bmatrix}
   F(I \otimes x^T) &  F(I \otimes u^T) &F \\ -F(I \otimes x^T) & -F(I \otimes u^T)& -F
\end{bmatrix} \\
Y \begin{bmatrix}
     \mathds{1}+\xi \\  \mathds{1}-\xi \\ \mathds{1} \\ \mathds{1}
\end{bmatrix}
&\leq  \lambda \psi(x) \mathds{1}.
\end{aligned}
\end{equation}
The proof follows now by replacing the constraint in \eqref{eq:conceptual} with the equivalent constraint \eqref{eq:duality}.
\end{pf}

\begin{remark} Note that \eqref{eq:duality} does not depend on the unknown dynamics $(A,B)$ or perturbation $v$. Thus, it avoids computing the vertices of $\mathcal{C}(\mathcal{D}_k)$, resulting in a substantial computational complexity reduction.
\end{remark}

\subsection{The Proposed DDRHC Algorithm}

Based on the results presented in the previous sections, we propose the  DDRHC algorithm  shown in Algorithm \ref{alg:main}.

\begin{algorithm}
\caption{Data-Driven Receding Horizon Control}\label{alg:main}
\begin{algorithmic}[1]
\Require Training data $\mathcal{D}_{train}$, constraints $\mathcal{X},\mathcal{U}, \mathcal{V}$
\State \textbf{Initialization} 
\begin{itemize}
\item [] Construct the consistency set $\mathcal{C}(\mathcal{D}_{train})$
\item []  Execute Algorithm \ref{alg:invariant} to compute $\mathcal{X}_I$, the largest controlled invariant set contained in $\mathcal{X}$
\item[] Set $k = 0,\; x_k = x_0$
\end{itemize}
\State \textbf{Execution} 
\Repeat
\State Solve \eqref{eq:dual_form} at $x_k$ to obtain $u_k$
\State Apply the control $u_k$ to the unknown system
\State Measure the resulting state $x_{k+1}$
\State Update $\mathcal{D}_{k+1} \gets \mathcal{D}_k \cup (x_k,u_k,x_{k+1})$
\State $k \gets k+1$
\Until{a stopping criterion is met}
\end{algorithmic}
\end{algorithm}

\begin{theorem}\label{teo:main}
Assume that the set $\mathcal{X}_I$ is controlled $\lambda$-contractive for some 
$0<\lambda \leq 1$ and for all systems in $\mathcal{C}(\mathcal{D}_{train})$. Then, 
the control action generated by Algorithm \ref{alg:main} renders the closed-loop system UUB in the set
$\lambda^*\mathcal{X}_I$ for some 
$0 \leq \lambda^* < \lambda$. Further, in the absence of disturbances, the origin is an exponentially stable equilibrium point, with a region of attraction $\mathcal{X}_I$.
\end{theorem}
\begin{pf} 
Follows from combining Lemma \ref{lemma:stable} and Theorem \ref{thm:conceptual}. 
\end{pf}

\section{Numerical Examples}
\label{sec:numerical_examples}

The following section presents the experimental results of the proposed data-driven receding horizon control method. All experiments were implemented in MATLAB 2024a with Yalmip \cite{Lofberg2004} and solved using Mosek \cite{mosek}. The supplementary code for the experiments is publicly available at 
\url{https://github.com/J-mzz/ddmpc}.

\begin{example}\label{eg:1}
    Consider a discrete-time linear system with
    \[
    A = \begin{bmatrix}
        0 & -0.99 \\ 0.99 & 0
    \end{bmatrix},\; 
    B = \begin{bmatrix}
        0 \\ 1
    \end{bmatrix}, 
    \]
    and state and input constraint $x\in\mathcal{X}=\left\{x \colon ||x||_\infty \leq 1 \right\}$, $u\in\mathcal{U}=\left\{u \colon ||u||_\infty \leq 1 \right\}$.  The system is subject to process noise, which is bounded by $\epsilon = 0.04$.
\end{example}

\begin{figure}[ht]
    \centering
    \begin{subfigure}[bt]{0.49\columnwidth}
        \centering
        \includegraphics[width=\linewidth]{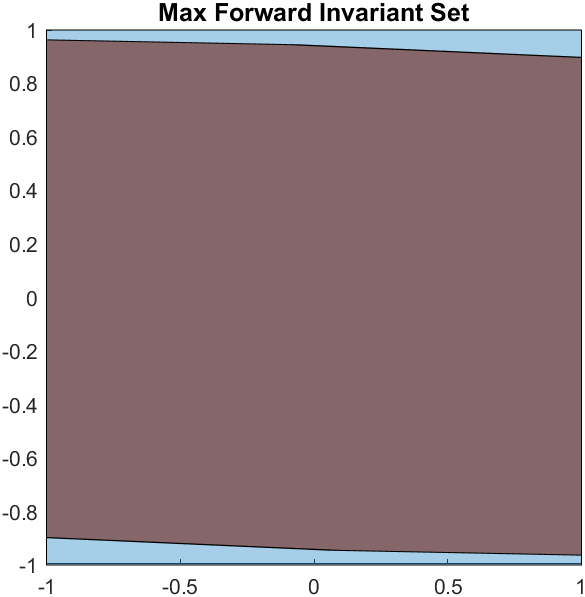}
        \caption{Largest Invariant Set}
        \label{fig:eg1_invariant}
    \end{subfigure}
    \hfill
    \begin{subfigure}[bt]{0.49\columnwidth}
        \centering
        \includegraphics[width=\linewidth]{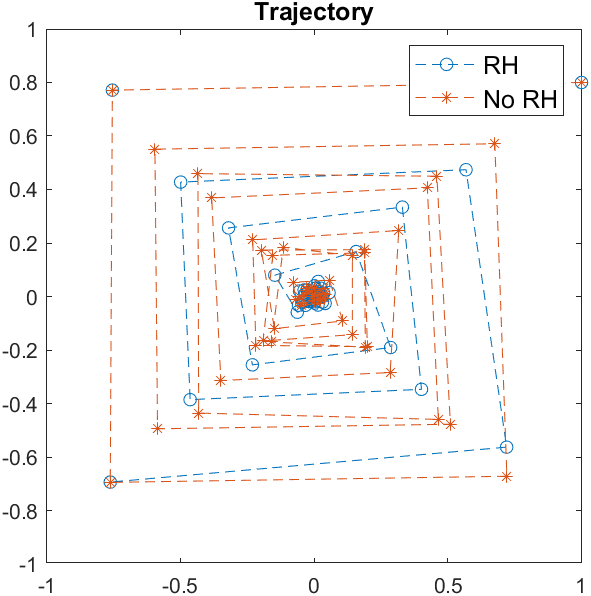}
        \caption{Trajectories with/without RH}
        \label{fig:eg1_traj}
    \end{subfigure}
    \vskip\baselineskip
    \begin{subfigure}[bt]{0.49\columnwidth}
        \centering
        \includegraphics[width=\linewidth]{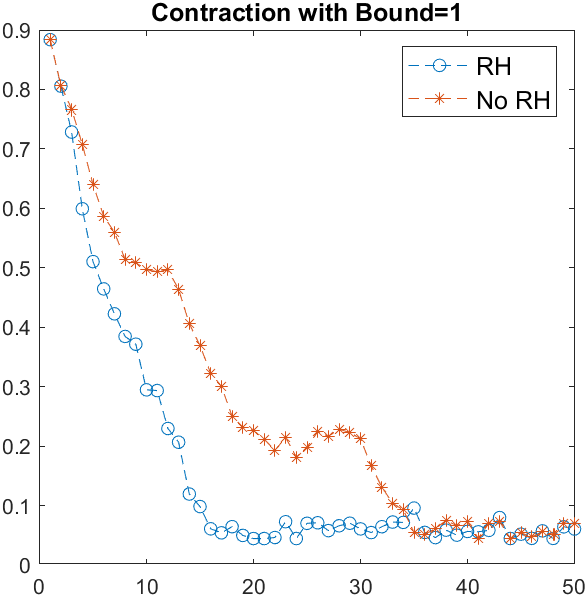}
        \caption{Contraction $\lambda$}
        \label{fig:eg1_lambda}
    \end{subfigure}
    \hfill
    \begin{subfigure}[bt]{0.49\columnwidth}
        \centering
        \includegraphics[width=\linewidth]{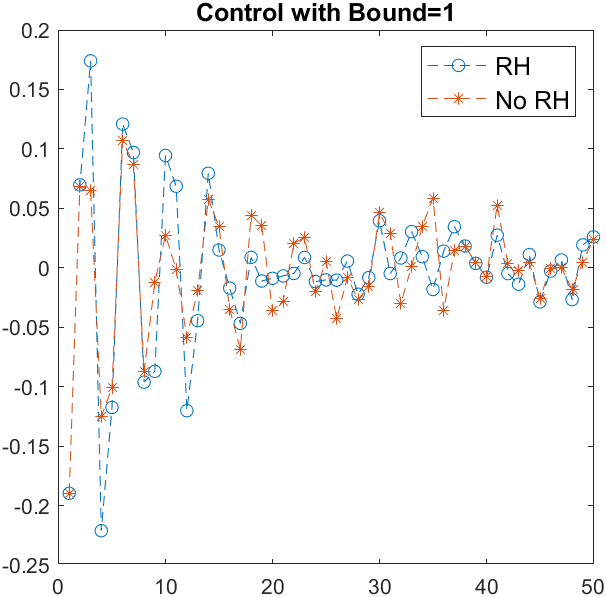}
        \caption{Control $u$}
        \label{fig:eg1_control}
    \end{subfigure}
    \caption{Results of Example 1}
    \label{fig:eg1}
\end{figure}

We collected 10 noisy measurements for the control design, and formulated a consistency set $\mathcal{C}(\mathcal{D}_0)$ from \eqref{eq:p1} of dimension 6, with 40 faces (of which 18 are nonredundant \cite{caron1989degenerate}) and 132 vertices. For the state constraint $\mathcal{X}$, which is the blue square in Fig. \ref{fig:eg1_invariant}, we applied Algorithm \ref{alg:invariant} for 20 iterations to obtain the 
largest controlled invariant set $\mathcal{X}_I\subseteq \mathcal{X}$, depicted as the brown polytope in Fig. \ref{fig:eg1_invariant}.

Starting from the initial point $x_0= [1; 0.8]$, which lies within the obtained invariant set, we applied our method for the unknown system and simulated it for 50 steps. Two trajectories, corresponding to our DDRHC and the regular DDC without RH, are presented in Fig. \ref{fig:eg1_traj}, while the corresponding $\lambda$ and $u$ are shown in Fig. \ref{fig:eg1_lambda} and \ref{fig:eg1_control}, respectively.

It is clear from Fig. \ref{fig:eg1_traj} that our method achieves better contraction performance compared to the regular DDC, as the blue line converges to the 
UUB set faster than the orange one. This is also shown  in Fig. \ref{fig:eg1_lambda} where the blue line remains below the orange line starting from the third step and until both lines reach the 
UUB set. This is because the data collected during the simulation further refine the consistency set and thus allow for more aggressive control, as illustrated in Fig. \ref{fig:eg1_control}.

\begin{example}\label{eg:2}
    Consider the same system as in Example \ref{eg:1}, with the same and state and input constraints $\mathcal{X}, \mathcal{U}$. But the noise level of the system is now increased to $\epsilon = 0.1$.
\end{example}

\begin{figure}[ht]
    \centering
    \begin{subfigure}[bt]{0.49\columnwidth}
        \centering
        \includegraphics[width=\linewidth]{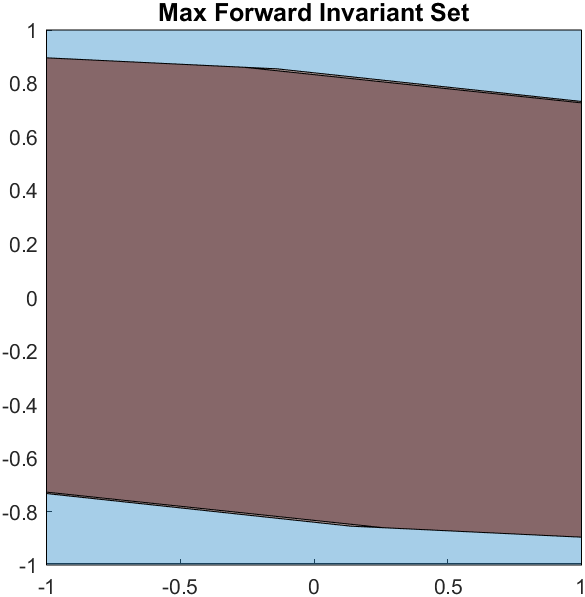}
        \caption{Largest Invariant Set}
        \label{fig:eg2_invariant}
    \end{subfigure}
    \hfill
    \begin{subfigure}[bt]{0.49\columnwidth}
        \centering
        \includegraphics[width=\linewidth]{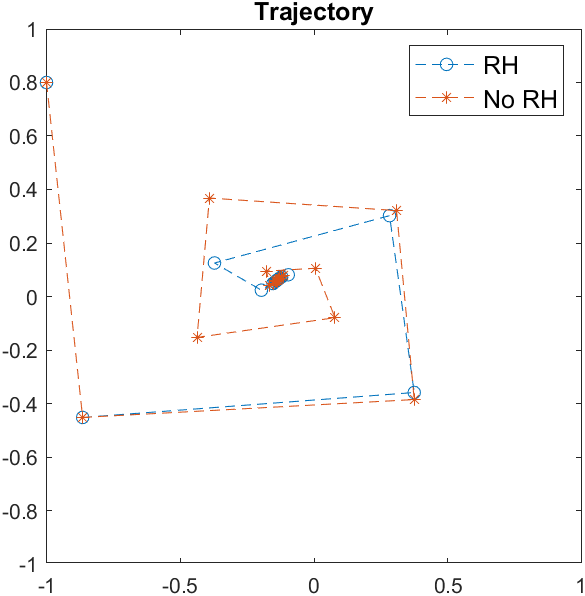}
        \caption{Trajectories with/without RH}
        \label{fig:eg2_traj}
    \end{subfigure}
    \vskip\baselineskip
    \begin{subfigure}[bt]{0.49\columnwidth}
        \centering
        \includegraphics[width=\linewidth]{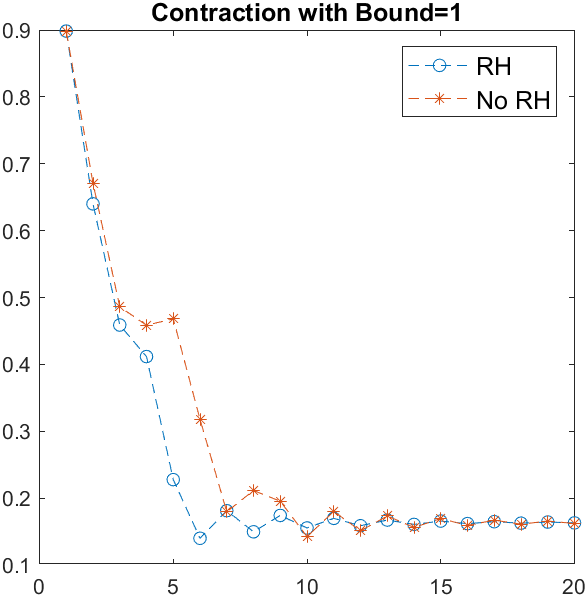}
        \caption{Contraction $\lambda$}
        \label{fig:eg2_lambda}
    \end{subfigure}
    \hfill
    \begin{subfigure}[bt]{0.49\columnwidth}
        \centering
        \includegraphics[width=\linewidth]{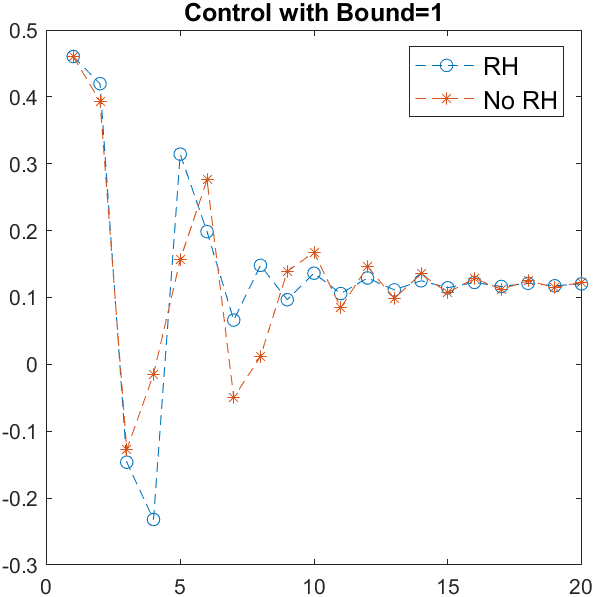}
        \caption{Control $u$}
        \label{fig:eg2_control}
    \end{subfigure}
    \caption{Results of Example 2}
    \label{fig:eg2}
\end{figure}

We collected 10 noisy measurements, and formulated a consistency set. For the state constraint $\mathcal{X}$, we applied Algorithm \ref{alg:invariant} for 20 iterations and obtained the 
largest controlled invariant set $\mathcal{X}_I\subseteq \mathcal{X}$, as shown in Fig. \ref{fig:eg2_invariant}.

Starting from the initial point $x_0= [-1; 0.8]$, we applied our method for the unknown system and simulated it for 20 steps. To better illustrate that our method achieves faster convergence, we applied a constant process noise in the simulation. Note that this is equivalent to analyzing a shifted system without process noise. Two trajectories, corresponding to our DDRHC and the regular DDC without RH, are presented in Fig. \ref{fig:eg2_traj}, while the corresponding $\lambda$ and $u$ are shown in Fig. \ref{fig:eg2_lambda} and \ref{fig:eg2_control}, respectively.

Fig. \ref{fig:eg2} clearly illustrates that the RH trajectory reaches the 
UUB set at time 6, while the non-RH  trajectory  does not reach the 
this set until time 9. 
We can also observe a smaller invariant set in Fig. \ref{fig:eg2_invariant}, compared to Fig. \ref{fig:eg1_invariant}, due to the increasing noise level.

\section{Conclusion}
\label{sec:conclusion}

In this work, we propose a novel data-driven receding horizon control framework for unknown discrete-time LTI systems, which guarantees robust performance in the presence of bounded process noise. Our method relies on set-membership constraints to achieve data-driven control and utilizes the online execution data to iteratively refine the set of all compatible systems. Experimental results demonstrate the effectiveness of the proposed framework and illustrate how receding horizon improves the regular data-driven control with faster convergence. Future work may focus on approximating the smallest controlled contractive set, optimizing the convergence over an extended prediction horizon, and relaxing the assumptions used in this work. 


\bibliography{reference.bib}             
                                                   







\appendix

\end{document}